\documentclass{article}
\usepackage{amsfonts}       
\usepackage{color} 
\newcommand{\email}[1]{{\small E-mail: {\textsf {#1}}}}
\newtheorem{theo}{Theorem}[section]
\newtheorem{prop}[theo]{Proposition}
\newtheorem{lem}[theo]{Lemma}

\newtheorem{rem}[theo]{Remark}

\newcommand{\mysection}[1]{\section{#1} \setcounter{equation}{0}}

\newcommand{\proof}{{\sc Proof.} \quad}

\newcommand{\be}{\begin{equation} \label}
\newcommand{\ee}{\end{equation}}
\newcommand{\bes}{\begin{equation} \begin{array}{c} \label}
\newcommand{\ees}{\end{array} \end{equation}}
\newcommand{\bea}{\begin{eqnarray}\label}
\newcommand{\eea}{\end{eqnarray}}
\newcommand{\beas}{\begin{eqnarray} \begin{array}{rcl} \label}
\newcommand{\eeas}{\end{array} \end{eqnarray}}
\newcommand{\bas}{\begin{eqnarray*}}\newcommand{\eas}{\end{eqnarray*}}
\newcommand{\bass}{\begin{eqnarray*} \begin{array}{rcl}}
\newcommand{\eass}{\end{array} \end{eqnarray*}}
\newcommand{\basss}{\begin{eqnarray*} \begin{array}{c}}
\newcommand{\easss}{\end{array} \end{eqnarray*}}
\newcommand{\qed}{{}\hfill //// \\[15pt]}
\newcommand{\bit}{\begin{itemize}}
\newcommand{\eit}{\end{itemize}}

\newcommand{\eps}{\varepsilon}
\newcommand{\abs}{\\[3mm]}

\newcommand{\dist}{{\rm dist}}
\newcommand{\divergence}{{\rm div}}
\newcommand{\trace}{{\rm trace}}
\begin{document}
\title{Convergence to steady states for radially symmetric solutions to a
quasilinear degenerate diffusive Hamilton-Jacobi equation\footnote{Partially supported by the ANR projects ``Hamilton-Jacobi et th\'eorie KAM faible'' (ANR-07-BLAN-3-187245) and EVOL (ANR-08-0242)}
}
\author{Guy Barles\footnote{Laboratoire de Math\'ematiques et Physique Th\'eorique, CNRS UMR~6083, F\'ed\'eration Denis Poisson, Universit\'e Fran\c{c}ois Rabelais, Parc de Grandmont, F--37200 Tours, France. \email{Guy.Barles@lmpt.univ-tours.fr}}, Philippe Lauren\c{c}ot\footnote{Institut de Math\'ematiques de Toulouse, CNRS UMR~5219, Universit\'e de Toulouse, F--31062 Toulouse cedex 9, France. \email{Philippe.Laurencot@math.univ-toulouse.fr}} \& 
Christian Stinner\footnote{Fachbereich Mathematik, Universit\"at Duisburg-Essen, D--45117 Essen, Germany. \email{Christian.Stinner@uni-due.de}}}
\date{\today}
\maketitle
\begin{abstract}
  \noindent Convergence to a single steady state is shown for non-negative and radially symmetric solutions to a diffusive Hamilton-Jacobi equation with homogeneous Dirichlet boundary conditions, the diffusion being the $p$-Laplacian operator, $p\ge 2$, and the source term a power of the norm of the gradient of $u$. As a first step, the radially symmetric and non-increasing stationary solutions are characterized. 
\end{abstract}
%
%
%
%
%
%
\mysection{Introduction}

We investigate the large time behaviour of non-negative and radially symmetric
solutions to the initial-boundary value problem 
\be{P}
        \left\{ \begin{array}{lllll}
        \partial_t u & = & \Delta_p u + |\nabla u|^q,   &     \qquad x\in B, & \ t\in (0,\infty), \\[1mm]
        u & = & 0,     &    \qquad x\in \partial B, & \ t\in (0,\infty),\\[2mm]
        u(x,0) & = & u_0(x),    &     \qquad x\in B, &
        \end{array} \right.
\ee
where $B:= \{ x \in \mathbb{R}^N \, : \, |x| < 1\}$ is the unit ball in
$\mathbb{R}^N$, $N \ge 2$, and the $p$-Laplacian operator is defined by
$$
\Delta_p u = \divergence (|\nabla u|^{p-2} \nabla u).
$$ 
We further assume the initial condition 
\begin{equation}\label{0.2}
  u_0 \in W^{1,\infty}_0 (B) \mbox{ is radially symmetric and non-negative and
} u_0 \not\equiv 0,
\end{equation} 
while the parameters $p$ and $q$ satisfy
\begin{equation}\label{0.3}
  p \ge 2 \quad\mbox{ and }\quad 0<q<p-1.
\end{equation}  

The partial differential equation in (\ref{P}) is a second-order parabolic equation
featuring a diffusion term (possibly quasilinear and degenerate if $p>2$) and a
source term $|\nabla u|^q$ counteracting the effect of diffusion and depending
solely on the gradient of the solution. The competition between the diffusion and
the source term is already revealed by the structure of steady states to (\ref{P}).
Indeed, while it follows from \cite[Theorem~1]{BB01} that zero is the only steady
state in $\mathcal{C}(\bar{B})$ when $p\ge 2$ and $q\ge p-1$, several steady states
may exist when $p\ge 2$ and $q\in (0,p-1)$ \cite{BDD92,Lau07,Sti09}. Another typical
feature of the competition between diffusion and source is the possibility of finite
time blow-up in a suitable norm, and this phenomenon has been shown to occur for
(\ref{P}) when $p=2$ and $q>2$ \cite{So02}. More precisely, it is established in
\cite{So02} that, when $p=2$ and $q>2$, there are classical solutions to (\ref{P})
  for which the $L^\infty$-norm of the gradient blows up in finite time, the
$L^\infty$-norm of the solution remaining bounded. These solutions may actually be
extended to all positive times in a unique way within the framework of viscosity
solutions \cite{BDL04,Tch09}, the boundary condition being also satisfied in the
viscosity sense. According to the latter, the homogeneous Dirichlet boundary
condition might not always be fulfilled for all times, a property which is likely to be
connected with the finite time blow-up of the gradient.

Coming back to the case where $p$ and $q$ fulfil (\ref{0.3}) and several steady
states may exist, a complete classification of steady states seems to be out of reach
when $B$ is replaced by an arbitrary open set of $\mathbb{R}^N$. Nevertheless, 
there are at least two situations in which the set of stationary solutions can be described,
namely, when $N=1$ and $B=(-1,1)$ \cite{Lau07,Sti09} and when $N\ge 2$ under the
additional requirement that the steady states are radially symmetric and
non-increasing, the latter being the first result of this paper. More precisely, we
show that (\ref{P}) has a one-parameter family of stationary solutions and that each
stationary solution is characterized by the value of its maximum.

\begin{theo}\label{theo1.1}
Assume (\ref{0.3}). Let $w \in W^{1,\infty} (B)$ be a radially symmetric and non-increasing viscosity
solution to $- \Delta_p w - |\nabla w|^q = 0$ in $B$ satisfying $w=0$ on $\partial B$. Then
there is $\vartheta \in [0,1]$ such that $w = w_\vartheta$, where 
  \begin{equation}\label{1.1.1}
    w_\vartheta (x) := c_0 \int\limits_{\max\{|x|, \vartheta\} }^1 \left(\rho -
\vartheta^\beta \rho^{-(\beta-1)} 
    \right)^{1 / (p-1-q)} {\rm d} \rho, \qquad \, x \in \bar{B},
  \end{equation}  
  for $\vartheta \in [0,1]$ with 
  \begin{equation}\label{1.1.2}
    \beta := 1+ \frac{(N-1)(p-1-q)}{p-1} >1 \quad\mbox{and}\quad c_0 := \left(
\frac{p-1-q}{(p-q)\beta} \right)^{1 /(p-1-q)}>0.
  \end{equation}  
  In particular, we have
  $w_0 (x) = (c_0 / \alpha)\  (1- |x|^\alpha )$ for $x \in \bar{B}$, where $\alpha
:= (p-q) /(p-1-q) >1$.
\end{theo}

\bigskip

\begin{rem}\label{rem1.1} As already mentioned, for any $M \in [0, c_0/\alpha]$
there is one and only one $\vartheta \in [0,1]$ such that $\| w_\vartheta
\|_{L^\infty (B)} = M$ as $\| w_\vartheta \|_{L^\infty (B)}$ is a decreasing
function of $\vartheta \in [0,1]$. This property plays an important role in the forthcoming analysis of the large time 
behaviour of solutions to (\ref{P}).
\end{rem}

\bigskip

Having a precise description of the set of steady states of (\ref{P}) at our
disposal, it is natural to investigate whether they attract the dynamics of
(\ref{P}) for large times. In other words, given a solution to (\ref{P}), does it
converge to a steady state as $t\to\infty$? A positive answer to this question is
given in \cite{Lau07,Sti09} when $N=1$, $B=(-1,1)$, and $p$ and $q$ fulfil
(\ref{0.3}). The one dimensional framework is fully exploited there as it allows the
construction of a Liapunov functional by the technique developed in \cite{Zel68}.
Such a nice tool does not seem to be available here and we instead use the theory of
viscosity solutions \cite{CIL92} and more precisely the relaxed half-limits method
introduced in \cite{BP88}. This approach has already been used in \cite{BS1,NR99,Ro01}
to investigate the large time behaviour of solutions to Hamilton-Jacobi equations and can be
roughly summarized as follows: given a non-negative and radially symmetric solution $u$ to
(\ref{P}) which is bounded in $W^{1,\infty}(B)$, the half-relaxed limits 
$$u_\ast (x) := \liminf\limits_{(s,\eps) \to (t,0)} u(x,\eps^{-1}s) \;\;\mbox{ and }\;\; u^\ast (x) := \limsup\limits_{(s,\eps) \to (t,0)} u(x,\eps^{-1}s), \quad x
\in \bar{B},$$
are well-defined, do not depend on $t > 0$, and are Lipschitz continuous viscosity supersolution and subsolution to 
$$
-\Delta_p z - |\nabla z|^q = 0 \;\;\;\mbox{ in }\;\;\; B , \quad z=0 \;\;\mbox{ on }\;\; \partial B ,
$$
respectively, by \cite[Lemma~6.1]{CIL92}. Clearly, $u_\ast\le u^\ast$ on $\bar{B}$ but we cannot apply the comparison principle at this stage to conclude that $u_\ast\ge u^\ast$ on $\bar{B}$. However, additional information are available  in this particular case, namely that $u_\ast$ and $u^\ast$ are both non-negative, radially symmetric, non-increasing, and have the same maximal value. Extensive use of these properties allows us to prove that $u_\ast\ge u^\ast$, from which we readily conclude that $u_\ast=u^\ast$ is a Lipschitz continuous radially symmetric and non-increasing stationary solution to (\ref{P}). Consequently, $u_\ast=u^\ast=w_\vartheta$ for some $\vartheta\in [0,1]$ by Theorem~\ref{theo1.1} and the assumption $u_0\not\equiv 0$ prevents $\vartheta=1$. The convergence result we obtain actually reads as follows.

\bigskip

\begin{theo}\label{theo1.2}
  Assume (\ref{0.2}) and (\ref{0.3}) and let $u$ denote the (radially symmetric)
viscosity solution to (\ref{P}). Then there
  is a unique $\vartheta \in [0,1)$ such that 
  $$\lim\limits_{t \to \infty} \| u(t) - w_\vartheta \|_{\mathcal{C} (\bar{B})} = 0.$$
\end{theo}
  
\bigskip

Notice that Theorem~\ref{theo1.2} applies in particular in the semilinear case $p=2$ with $q\in (0,1)$ according to (\ref{0.3}). Still in the semilinear case $p=2$, several results on the large time behaviour of solutions to (\ref{P}) are also available when $q\ge 1$ and $B$ is replaced by an arbitrary open set $\Omega$ of $\mathbb{R}^N$ \cite{ARBS04,BDHL07,SZ06,Tch09}, including the convergence to zero of global solutions which are bounded in $W^{1,\infty}(\Omega)$.

\bigskip

The analysis in this paper being restricted to radially symmetric solutions, we define $r := |x|$ and switch between
the notation $u=u(x,t)$ and $u=u(r,t)$, whenever this is convenient.

For further use, we introduce the following notations:
\begin{equation}\label{1.6}
  F(s,X) := - |s|^{p-2} \trace (X) - (p-2) |s|^{p-4} \langle X s,s \rangle - |s|^q
\quad\mbox{ for } (s,X) \in \mathbb{R}^N\times \mathbb{R}^{N \times N}  ,
\end{equation}
its radially symmetric counterpart 
\begin{equation}\label{1.7}
  f(r,\mu, \zeta) := - (p-1) |\mu|^{p-2} \zeta - \frac{N-1}{r} |\mu|^{p-2} \mu -
|\mu|^q \quad\mbox{for } (r,\mu,\zeta) \in (0,1)\times\mathbb{R}\times\mathbb{R}, 
\end{equation}
and the radially symmetric $p$-Laplacian operator
\begin{equation}\label{1.7b}
  f_0(r,\mu, \zeta) := - (p-1) |\mu|^{p-2} \zeta - \frac{N-1}{r} |\mu|^{p-2} \mu
 \quad\mbox{for } (r,\mu,\zeta) \in (0,1)\times\mathbb{R}\times\mathbb{R}. 
\end{equation}
%
%
\mysection{Radially symmetric and non-increasing stationary solutions}

In this section, we prove Theorem~\ref{theo1.1}, that is, if $w$ is a radially
symmetric, non-increasing, and Lipschitz
continuous viscosity solution to the stationary equation
\begin{equation}\label{2.1}
  \left\{ \begin{array}{rl} - \Delta_p w - |\nabla w|^q = 0 &\qquad\mbox{in } B, \\
  w = 0 & \qquad\mbox{on } \partial B, \end{array} \right.
\end{equation}  
then $w = w_\vartheta$ for some $\vartheta \in [0,1]$. To this end, we first observe
that, as a function of $r = |x|$,
$w$ is a viscosity solution to $f(r, \partial_r w, \partial_r^2 w) = 0$ in $(0,1)$
with $w(1) = 0$ (recall that $f$ is
defined in (\ref{1.7})). 

Next, as a preliminary step, let us first give a formal proof, assuming $w$ to be in
$\mathcal{C}^1 (\bar{B})$ and solving
(\ref{2.1}) pointwise. In particular, we will derive an identity (see (\ref{2.3})
below) which turns out to be valid for
viscosity solutions as we shall see later on. 

As $w$ is radially symmetric and in $\mathcal{C}^1 (\bar{B})$, we have $\partial_r w(0) = 0$.
In addition, by (\ref{2.1}), 
  $$\varphi (r) := r^{N-1} (|\partial_r w|^{p-2} \partial_r w)(r), \quad r \in [0,1],$$
  fulfils $\varphi \in W^{1,\infty} ((0,1))$ with $\partial_r \varphi (r) = -
r^{N-1}|\partial_r w(r)|^q \le 0$ a.e. in $(0,1)$.
  Thus, $\varphi$ is a non-increasing function in $[0,1]$. As moreover $w$ is
non-increasing with $w(1)=0$, we have $\partial_r w(1) \le 0$.
  
  Now, either $\partial_r w (1) =0$ and thus $\varphi (1)=0$. Since $\varphi$ is
non-increasing with $\varphi(0) = 0$, we conclude that $\varphi \equiv 0$. This
implies $w = w_1 \equiv 0$. 
  
  Or $\partial_r w (1) <0$, and the continuity and monotonicity of $\varphi$
warrant that there is a unique $\vartheta \in [0,1)$ such that
  $\varphi =0$ in $[0, \vartheta]$ and $\varphi <0$ in $(\vartheta,1]$.
  Hence, 
  $$
  \partial_r \varphi(r) = - r^{[(N-1)(p-1-q)] / (p-1)} |\varphi(r)|^{q / (p-1)} = -
r^{\beta -1}
  (-\varphi(r))^{q / (p-1)} \;\;\mbox{ in }\;\; (\vartheta,1).
  $$
  After integration we obtain
  $$-\frac{p-1}{p-1-q} (-\varphi (r))^{(p-1-q) / (p-1)} + \frac{1}{\beta} r^\beta =
\gamma \quad\mbox{for } r \in
  (\vartheta,1)$$
  with some constant $\gamma \in \mathbb{R}$. Introducing
  \begin{equation}\label{2.2}
    \chi(z) := \frac{p-1}{p-1-q} |z|^{p-2-q} z \qquad\mbox{for } z \in \mathbb{R},
  \end{equation}  
  we end up with
  \begin{equation}\label{2.3}
    r^{\beta-1} \chi(\partial_r w(r)) + \frac{1}{\beta} r^\beta = \gamma \quad\mbox{for
} r \in (\vartheta,1)
  \end{equation}
  as $\partial_r w <0$ in $(\vartheta,1)$. Letting $r \searrow \vartheta$ implies
$\gamma = \vartheta^\beta/\beta$ owing to $\partial_r w(\vartheta) =0$ and
$0<q<p-1$. 
  
  Furthermore, due to $\partial_r w <0$ in $(\vartheta,1)$, we have
  $$-\frac{p-1}{p-1-q} \left( r^{(N-1) / (p-1)} (-\partial_r w (r)) \right)^{p-1-q}
= \frac{1}{\beta} \left(\vartheta^\beta -
  r^\beta \right) \quad\mbox{for } r \in (\vartheta,1).$$
  Hence, we conclude
  $$\partial_r w (r) = - \left( \frac{p-1-q}{(p-1) \beta} \left( r- \vartheta^\beta
r^{-(\beta-1)} \right) 
    \right)^{1 / (p-1-q)} \quad\mbox{for } r \in (\vartheta,1).$$
  Using $w(1) =0$ and the definition of $c_0$, a further integration implies
  $$w (r) = c_0 \int\limits_r^1 \left(\rho - \vartheta^\beta \rho^{-(\beta-1)} 
  \right)^{1 / (p-1-q)} {\rm d} \rho = w_\vartheta(r) \qquad\mbox{for } r \in
[\vartheta,1].$$
  Furthermore, we get $w(r) = w(\vartheta)$ for any $r \in [0,\vartheta]$ since
$\partial_r w \equiv 0$ in $[0,\vartheta]$ and we conclude that $w = w_\vartheta$.

\bigskip

We now turn to the proof of Theorem~\ref{theo1.1} and first establish some
preliminary results. We recall that, by the Rademacher theorem, a Lipschitz
continuous function $v\in W^{1,\infty}((0,1))$ is differentiable a.e. and the measure
of the differentiability set
$$D(v) := \{ r_0 \in (0,1) \, : \, \partial_r v(r_0) \mbox{ exists } \}$$
is thus equal to one.

\begin{lem}\label{lem2.1}
  Let $v \in W^{1,\infty} ((0,1))$ be a non-negative and non-increasing viscosity
supersolution to
  \begin{equation}
  \label{2.3b}
  f_0(r,\partial_r z, \partial_r^2 z) = 0 \quad \mbox{in } (0,1),
  \end{equation}
  the Hamiltonian $f_0$ being defined in (\ref{1.7b}). Then, if $r_1\in D(v)$ and
$r_2\in D(v)$ are such that $r_1<r_2$, we have
  $$r_2^{(N-1) / (p-1)} \partial_r v(r_2) \le r_1^{(N-1) / (p-1)} \partial_r v(r_1) .$$
\end{lem}

\proof Take $0<r_1 < r_2 <1$ with $r_1,r_2 \in D(v)$ and assume for contradiction that
  $$\xi_1 := r_1^{(N-1) / (p-1)} \partial_r v(r_1) < r_2^{ (N-1) / (p-1)} \partial_r
v(r_2) =: \xi_2.$$
  As $v$ is non-increasing we have $\xi_2 \le 0$. Now take $\xi_1 < \eta_1 < \eta_2
< \xi_2 \le 0$ and define 
  $\Phi$ by
  $$r^{(N-1) / (p-1)} \partial_r \Phi (r) = \eta_1 + (\eta_2 - \eta_1)
\frac{r-r_1}{r_2-r_1}, \qquad r \in [r_1,r_2],$$
  along with $\Phi (r_1) = 0$. 

  On the one hand, $v - \Phi$ is continuous in $[r_1,r_2]$ and thus attains its
minimum at a point $r_0 \in
  [r_1,r_2]$. On the other hand, we have
  $$\partial_r (v- \Phi) (r_1) = \frac{\xi_1 - \eta_1}{r_1^{(N-1) / (p-1)}} <0
\quad\mbox{and}\quad
    \partial_r (v- \Phi) (r_2) = \frac{\xi_2 - \eta_2}{r_2^{(N-1) / (p-1)}} >0$$
  so that we cannot have $r_0 =r_1$ or $r_0 = r_2$. Thus, $r_0 \in (r_1,r_2)$ and,
since $v$ is a viscosity 
  supersolution to (\ref{2.3b}), we have 
  $$ -\frac{1}{r_0^{N-1}} \partial_r \left( r^{N-1} |\partial_r \Phi|^{p-2} 
    \partial_r \Phi \right) (r_0) \ge 0.$$
  Since $r^{(N-1) / (p-1)}  \partial_r \Phi(r) \le \eta_2<0$ for $r\in [r_1,r_2]$ we obtain
  \begin{eqnarray*}
    - \left( r^{N-1} |\partial_r \Phi|^{p-2} \partial_r \Phi \right)(r) &=& r^{N-1}
|\partial_r \Phi (r)|^{p-1} 
    = \left( - r^{(N-1)/(p-1)} \partial_r \Phi(r) \right)^{p-1} \\ 
    &=& \left| \eta_1 + (\eta_2 - \eta_1) \frac{r-r_1}{r_2-r_1} \right|^{p-1}.
  \end{eqnarray*}  
  Differentiating and taking $r=r_0$, we end up with
  \begin{eqnarray*}
    0 &\le& - \partial_r \left( r^{N-1} |\partial_r \Phi|^{p-2} \partial_r \Phi
\right) (r_0) \\
    &=& (p-1) \left| \eta_1 + (\eta_2 - \eta_1) \frac{r_0-r_1}{r_2-r_1}
\right|^{p-3} \left( \eta_1 + (\eta_2 - \eta_1) 
    \frac{r_0-r_1}{r_2-r_1} \right) \ \frac{\eta_2 - \eta_1}{r_2-r_1} < 0,
  \end{eqnarray*}  
  and a contradiction. \qed

\medskip

In order to show that a viscosity solution to (\ref{2.1}) satisfies (\ref{2.3}), we
next prove that the
left-hand side of (\ref{2.3}) is non-increasing for a supersolution to (\ref{2.1}). 

\begin{lem}\label{lem2.2}
  Let $w \in W^{1,\infty} ((0,1))$ be a non-increasing viscosity supersolution to
  $f(r, \partial_r z, \partial_r^2 z) = 0$ in $(0,1)$ such that $\| w \|_{L^\infty
  ((0,1))} >0$ and $w(1)=0$, and define $r_0\in [0,1]$ by
  $$
  r_0 := \inf\left\{r \in (0,1] \, : \, w(r) < \|w\|_{L^\infty ((0,1))} \right\}.
  $$ 
  If $r_1\in D(w)$ and $r_2\in D(w)$ are such that $r_0<r_1<r_2$, then 
  $$
  r_1^{\beta -1} \chi(\partial_r w(r_1)) + \frac{r_1^\beta}{\beta} \ge 
  r_2^{\beta -1} \chi(\partial_r w(r_2)) + \frac{r_2^\beta}{\beta},
  $$
  the parameter $\beta$ and the function $\chi$ being defined in (\ref{1.1.2}) and
(\ref{2.2}), respectively.
\end{lem}

\proof The properties of $w$ imply $r_0 \in [0,1)$. As $w$ is non-increasing and
Lipschitz continuous, the definition of $r_0$ yields that there is a sequence
$(\varrho_n)_{n \ge 1}$ such that $\varrho_n \in D(w)$, $\partial_r w(\varrho_n) <0$
and $\varrho_n \searrow r_0$ as $n \to \infty$. Pick $r_1\in D(w)\cap (r_0,1)$. For
$n$ large enough, we have $r_1>\varrho_n$. Since $w$ is clearly also a supersolution
to (\ref{2.3b}), we infer from Lemma~\ref{lem2.1} that 
$$
r_1^{(N-1) / (p-1)} \partial_r w(r_1) \le \varrho_n^{(N-1) / (p-1)} \partial_r
w(\varrho_n)< 0
$$ 
for $n$ large enough. Consequently,
  \begin{equation}\label{2.2.1}
  r_1^{(N-1) / (p-1)} \partial_r w(r_1) <0 \quad \mbox{for } r_1\in D(w)\cap (r_0,1).
  \end{equation}  

Assume now for contradiction that there are $r_1, r_2 \in (r_0,1)
  \cap D(w)$ such that $r_1 <r_2$ and
  $$r_1^{\beta -1} \chi(\partial_r w(r_1)) + \frac{r_1^\beta}{\beta} < 
  r_2^{\beta -1} \chi(\partial_r w(r_2)) + \frac{r_2^\beta}{\beta}.$$
  As $\partial_r w(r_1) <0$ by (\ref{2.2.1}), we have $\chi(\partial_r w(r_1))<0$ and
we can choose two real numbers $\eta_1$ and $\eta_2$ such that
  $$r_1^{\beta -1} \chi(\partial_r w(r_1)) + \frac{r_1^\beta}{\beta} < \eta_1 < \eta_2 < 
  r_2^{\beta -1} \chi(\partial_r w(r_2)) + \frac{r_2^\beta}{\beta}\ , \qquad \eta_1 <
\frac{r_1^\beta}{\beta}, $$
  and 
  $$a:= 1 - \frac{\beta (\eta_2 - \eta_1)}{r_2^\beta - r_1^\beta} \in (0,1).$$
Indeed we first choose $ \eta_1 \in (r_1^{\beta -1} \chi(\partial_r w(r_1)) + (r_1^\beta/\beta), r_1^\beta/\beta)$ and then $\eta_2>\eta_1$ close enough to $\eta_1$ in order to have $a\in (0,1)$.  Setting now
  $$A := \eta_1 - (1-a) \frac{r_1^\beta}{\beta} = \eta_2 - (1-a)
\frac{r_2^\beta}{\beta},$$
  let $\Phi$ denote the solution to
  \begin{equation}\label{2.2.2}
    r^{\beta-1} \chi(\partial_r \Phi(r)) + a \frac{r^\beta}{\beta} = A, \qquad r \in
[r_1,r_2],    
  \end{equation} 
  such that $\Phi (r_1) =0$. Observe that the choice of $a$ and $A$ ensure that  and
  \begin{equation}\label{2.2.3}
    r_i^{\beta-1} \chi(\partial_r \Phi(r_i)) + \frac{r_i^\beta}{\beta} = \eta_i
\qquad\mbox{for } i=1,2.
  \end{equation}
  Due to 
  $$A - a \frac{r_1^\beta}{\beta} = \eta_1 - \frac{r_1^\beta}{\beta} <0$$
  we conclude by (\ref{2.2.2}) that
  $$\chi(\partial_r \Phi(r)) = r^{-(\beta-1)} \left( A - a \frac{r^\beta}{\beta}
\right) \le r^{-(\beta-1)} \left( A - a \frac{r_1^\beta}{\beta} \right) <0
\qquad\mbox{for } r \in [r_1,r_2].$$
  This implies that $\partial_r \Phi(r) <0$ for $r\in [r_1,r_2]$, so that $\Phi \in
\mathcal{C}^2([r_1,r_2])$ by (\ref{2.2.2}). In addition,
  $$ \quad (- \partial_r \Phi(r))^{p-1-q} = \frac{p-1-q}{p-1} \left( \frac{a}{\beta} r
  - A r^{-(\beta -1)} \right), \quad r \in [r_1,r_2],$$
  hence
  $$\partial_r \Phi(r) = - \left[ \frac{p-1-q}{p-1} \left( \frac{a}{\beta} r - A
r^{-(\beta -1)} \right) \right]^{1 /
    (p-1-q)}, \qquad r \in [r_1, r_2].$$
  Furthermore, due to (\ref{2.2.3}) and the choice of $\eta_1$, we obtain
  $$r_1^{\beta-1} \chi(\partial_r w(r_1)) + \frac{r_1^\beta}{\beta} < \eta_1 =
r_1^{\beta-1} \chi(\partial_r \Phi(r_1)) 
    + \frac{r_1^\beta}{\beta}.$$
  This implies $\chi(\partial_r w(r_1)) < \chi(\partial_r \Phi(r_1))$ and, since $\chi$ is
increasing,
  $$\partial_r w(r_1) < \partial_r \Phi(r_1).$$
  Similarly, we conclude
  $$\partial_r w(r_2) > \partial_r \Phi(r_2).$$  
  Now $w - \Phi$ is a continuous function in $[r_1, r_2]$ and thus attains its minimum
  at some $r_m \in [r_1,r_2]$. The above two inequalities prevent $r_m$ to be equal
to $r_1$ or $r_2$
  and, since $w$ is a viscosity supersolution to $f(r, \partial_r v, \partial_r^2 v)
= 0$ in $(0,1)$, we have
  $$-\frac{1}{r_m^{N-1}} \partial_r \left( r^{N-1} |\partial_r \Phi (r)|^{p-2}
\partial_r \Phi(r) \right)(r_m) 
    - | \partial_r \Phi(r_m)|^q \ge 0.$$ 
  But as $\partial_r \Phi <0$, (\ref{2.2.2}) implies
  \begin{eqnarray}\label{2.2.4}
    & & - \partial_r \left( r^{N-1} |\partial_r \Phi (r)|^{p-2} \partial_r \Phi(r)
\right) =
    \partial_r \left( r^{N-1} |\partial_r \Phi (r)|^{p-1} \right) \nonumber \\
    &=& \partial_r \left( \left| \frac{p-1-q}{p-1} r^{\beta-1} \chi(\partial_r \Phi(r))
\right|^{(p-1) / (p-1-q)} \right) 
    \nonumber \\
    &=& - a r^{\beta-1} \left| \frac{p-1-q}{p-1} r^{\beta-1} \chi(\partial_r \Phi(r))
\right|^{[(p-1) / (p-1-q)] -2}
    \left( \frac{p-1-q}{p-1} r^{\beta-1} \chi(\partial_r \Phi(r)) \right) \nonumber \\
    &=& a r^{(\beta-1)(p-1)/(p-1-q)} \left| \frac{p-1-q}{p-1} \chi(\partial_r \Phi(r))
\right|^{[(p-1) / (p-1-q)]-1}
    \nonumber \\
    &=& a r^{N-1} |\partial_r \Phi(r)|^q \qquad\mbox{for } r \in [r_1,r_2],
  \end{eqnarray}  
  so that 
  $$-\frac{1}{r_m^{N-1}} \partial_r \left( r^{N-1} |\partial_r \Phi (r)|^{p-2}
\partial_r \Phi(r) \right)(r_m) 
    - | \partial_r \Phi(r_m)|^q = (a-1) |\partial_r \Phi (r_m)|^q < 0$$
  since $a<1$, and a contradiction. \qed

In a similar way we now establish that the left-hand side of (\ref{2.3}) is
non-decreasing for viscosity subsolutions to 
(\ref{2.1}). 

\begin{lem}\label{lem2.3}
  Let $w \in W^{1,\infty} ((0,1))$ be a non-increasing viscosity subsolution to
  $f(r, \partial_r z, \partial_r^2 z) = 0$ in $(0,1)$ such that $\| w \|_{L^\infty
  ((0,1))} >0$ and $w(1)=0$, and define $r_0\in [0,1]$ by
  $$
  r_0 := \inf\left\{r \in (0,1] \, : \, w(r) < \|w\|_{L^\infty ((0,1))} \right\}.
  $$ 
  If $r_1\in D(w)$ and $r_2\in D(w)$ are such that $r_0<r_1<r_2$, then 
  $$r_1^{\beta -1} \chi(\partial_r w(r_1)) + \frac{r_1^\beta}{\beta} \le 
  r_2^{\beta -1} \chi(\partial_r w(r_2)) + \frac{r_2^\beta}{\beta}.$$
\end{lem}

\proof The properties of $w$ imply $r_0 \in [0,1)$. Assume for contradiction that
there are $r_1, r_2 \in (r_0,1)
  \cap D(w)$ such that $r_1 <r_2$ and
  $$r_1^{\beta -1} \chi(\partial_r w(r_1)) + \frac{r_1^\beta}{\beta} > 
  r_2^{\beta -1} \chi(\partial_r w(r_2)) + \frac{r_2^\beta}{\beta}.$$
  We may then choose $\eta_1, \eta_2 \in \mathbb{R}$
  such that
  $$r_1^{\beta -1} \chi(\partial_r w(r_1)) + \frac{r_1^\beta}{\beta} > \eta_1 > \eta_2 > 
  r_2^{\beta -1} \chi(\partial_r w(r_2)) + \frac{r_2^\beta}{\beta},$$
  and define
  $$a:= 1 + \frac{\beta (\eta_1 - \eta_2)}{r_2^\beta - r_1^\beta} >1 \quad \mbox{
and }\quad A := \eta_1 + (a-1) \frac{r_1^\beta}{\beta} = \eta_2 + (a-1)
\frac{r_2^\beta}{\beta}.$$
  Let $\Phi$ denote the solution to
  \begin{equation}\label{2.3.2}
    r^{\beta-1} \chi(\partial_r \Phi(r)) + a \frac{r^\beta}{\beta} = A, \qquad r \in
[r_1,r_2],    
  \end{equation} 
  such that $\Phi (r_1) =0$. Thanks to the choice of $a$ and $A$, we have
  \begin{equation}\label{2.3.3}
    r_i^{\beta-1} \chi(\partial_r \Phi(r_i)) + \frac{r_i^\beta}{\beta} = \eta_i
\qquad\mbox{for } i=1,2,
  \end{equation}
  and the monotonicity of $w$ implies that
  $$A - a \frac{r_1^\beta}{\beta} = \eta_1 - \frac{r_1^\beta}{\beta} < r_1^{\beta-1}
\chi(\partial_r w(r_1)) \le 0.$$
  Consequently,
  $$\chi(\partial_r \Phi(r)) = r^{-(\beta-1)} \left( A - a \frac{r^\beta}{\beta}
\right) \le r^{-(\beta-1)} \left( A - a \frac{r_1^\beta}{\beta} \right) <0
\qquad\mbox{for } r \in [r_1,r_2],$$
  hence $\partial_r \Phi(r) <0$ for $r \in [r_1,r_2]$. We then conclude from
(\ref{2.3.2}) that $\Phi \in \mathcal{C}^2 ([r_1,r_2])$. 
  Furthermore, due to (\ref{2.3.3}), the choices of $\eta_1$ and $\eta_2$, and the
monotonicity of $\chi$, we obtain
  $$\partial_r w(r_1) > \partial_r \Phi(r_1) \qquad\mbox{and}\qquad \partial_r
w(r_2) < \partial_r \Phi(r_2).$$
  
  Now $w - \Phi$ is a continuous function in $[r_1,r_2]$ and thus attains its maximum
  at some point $r_m \in [r_1,r_2]$. The above two inequalities prevent $r_m$ to be
equal to $r_1$ or $r_2$
  and, since $w$ is a viscosity subsolution to $f(r, \partial_r v, \partial_r^2 v) =
0$ in $(0,1)$, we have
  $$-\frac{1}{r_m^{N-1}} \partial_r \left( r^{N-1} |\partial_r \Phi (r)|^{p-2}
\partial_r \Phi(r) \right)(r_m) 
    - | \partial_r \Phi(r_m)|^q \le 0.$$ 
  But, owing to $\partial_r \Phi(r) <0$, (\ref{2.3.2}) and $a>1$,   
  we conclude similarly to (\ref{2.2.4}) that
  $$-\frac{1}{r_m^{N-1}} \partial_r \left( r^{N-1} |\partial_r \Phi (r)|^{p-2}
\partial_r \Phi(r) \right)(r_m) 
    - | \partial_r \Phi(r_m)|^q = (a-1) |\partial_r \Phi (r_m)|^q > 0$$
  and end up with a contradiction. 
\qed

We are now in a position to prove Theorem~\ref{theo1.1}. The keystone of the proof
is that, according to Lemma~\ref{lem2.2} and Lemma~\ref{2.3}, any non-increasing
viscosity solution to $f(r, \partial_r v, \partial_r^2 v) = 0$ in $(0,1)$ satisfying
$w(1)=0$ has to fulfil (\ref{2.3}).

\bigskip

{\sc Proof of Theorem~\ref{theo1.1}.}

  Let $w \in W^{1,\infty} ((0,1))$ be a non-increasing viscosity solution to $f(r,
\partial_r v, \partial_r^2 v) = 0$ in $(0,1)$ satisfying $w(1)=0$. Either $w
\equiv 0=w_1$ or $M:= \| w \|_{L^\infty ((0,1))} >0$ and we define $r_0 \in [0,1)$ by
$$r_0 := \inf\{r \in (0,1] \, : \, w(r) < M \}.$$ 

  Now, owing to Lemma~\ref{lem2.2} and Lemma~\ref{lem2.3}, there is a constant
$\gamma \in \mathbb{R}$ such that 
  \begin{equation}\label{t1.1.1}
    r^{\beta -1} \chi(\partial_r w(r)) + \frac{r^\beta}{\beta} = \gamma
  \end{equation}  
  for any $r \in (r_0,1) \cap D(w)$ and thus a.e. in $(r_0,1)$. Combining the
monotonicity of $w$ and $\chi$ with (\ref{t1.1.1}), we moreover deduce that  
  \begin{equation}\label{t1.1.2}
    \gamma \le \frac{r_0^\beta}{\beta}
  \end{equation}  
  and
  $$\partial_r w(r) = - \left[ \frac{p-1-q}{p-1} \left( \frac{r}{\beta} - \gamma
r^{-(\beta -1)} \right) 
    \right]^{1 / (p-1-q)} \qquad\mbox{for a.e. } r \in (r_0,1).$$
  Integrating and using the boundary condition $w(1)=0$, we obtain
  $$w(r) =  \int\limits_r^1 \left[ \frac{p-1-q}{(p-1)\beta} \left(\rho - \gamma
\beta \rho^{-(\beta-1)} 
    \right) \right]^{1 / (p-1-q)} {\rm d} \rho \qquad\mbox{for any } r \in [r_0,1].$$
  Recalling $w(r) \equiv M$ for $r \in [0,r_0]$ and the definition of $c_0$, we
conclude
  \begin{equation}\label{t1.1.3}
    w(r) =  c_0 \int\limits_{\max\{r, r_0\} }^1 \left(\rho - \gamma \beta
\rho^{-(\beta-1)} 
    \right)^{1 / (p-1-q)} {\rm d} \rho, \qquad r \in [0,1].
  \end{equation}
  
  It remains to show that $\gamma = r_0^\beta/\beta$ in order to obtain that $w
= w_{r_0}$.\\
  Consider first the case $r_0 =0$. Since $\beta>1$, the Lipschitz continuity of
$w$ yields $\gamma = 0 =   r_0^\beta/\beta$ by letting $r \searrow 0$ in
(\ref{t1.1.1}). \\
  Next, if $r_0 \in (0,1)$, we assume for contradiction that $\gamma <
r_0^\beta/\beta$. Then we fix
  $\vartheta \in [0,r_0)$ such that $\gamma < \vartheta^\beta/\beta$ and
choose $\Lambda >1$ such that
  $$\Lambda^{p-1-q} < 1 + \vartheta^\beta - \gamma \beta.$$
  This choice of $\Lambda$ implies that the function
  $$g(r) := \left( 1 - \gamma \beta r^{-\beta} \right) - \Lambda^{p-1-q} \left( 1-
\vartheta^\beta r^{-\beta} \right),
  \quad r \in (r_0,1),$$
  satisfies
  $$g^\prime (r) =  \beta^2 r^{-\beta-1} \left( \gamma - \Lambda^{p-1-q}
\frac{\vartheta^\beta}{\beta} \right)
    \le \beta^2 r^{-\beta-1} \left( \gamma - \frac{\vartheta^\beta}{\beta} \right)
<0, \quad r \in (r_0,1),$$
  and thus
  $$g(r) \ge g(1) \ge 1 - \gamma \beta - \Lambda^{p-1-q} + \vartheta^\beta >0, \quad
r\in [r_0,1].$$  
  Consequently,
  $$
  \left( 1 - \gamma \beta r^{-\beta} \right) > \Lambda^{p-1-q} \left( 1-
\vartheta^\beta r^{-\beta} \right),
  \quad r \in [r_0,1],
  $$
  and it follows from (\ref{t1.1.3}) that
  \begin{eqnarray*}
    \partial_r w (r) &=& - c_0 r^{1/(p-1-q)} \left(1 - \gamma \beta r^{-\beta}
\right)^{1/ (p-1-q)} \\ &<& 
    - c_0 r^{1/(p-1-q)} \Lambda \left( 1- \vartheta^\beta r^{-\beta} \right)^{1 /
(p-1-q)} = \Lambda \partial_r w_\vartheta (r), \quad r \in (r_0,1).
  \end{eqnarray*}  
  In particular, $w(r)-\Lambda w_\vartheta(r) \le w(r_0)-\Lambda w_\vartheta(r_0)$
for $r\in [r_0,1]$. Furthermore, 
  $$w (r) - \Lambda w_\vartheta (r) = w(r_0) - \Lambda w_\vartheta (r) \le w(r_0) -
\Lambda w_\vartheta (r_0),
    \quad r \in [0, r_0],$$
thanks to the monotonicity of $w_\vartheta$, and the function $w -
\Lambda w_\vartheta$ has a global maximum at $r_0$. Since $w_\vartheta \in \mathcal{C}^2 (
(\vartheta,1))$, $\vartheta < r_0$, and $w$ is a viscosity subsolution to $f(r,
\partial_r v, \partial_r^2 v) = 0$ in $(0,1)$, we conclude that
  $$f(r_0, \partial_r (\Lambda w_\vartheta) (r_0), \partial_r^2 (\Lambda
w_\vartheta)(r_0)) \le 0.$$ 
  However, as $\Lambda >1$ and $\vartheta < r_0$, we clearly have
  $$f(r_0, \partial_r (\Lambda w_\vartheta) (r_0), \partial_r^2 (\Lambda
w_\vartheta)(r_0))
    = \left(\Lambda^{p-1} - \Lambda^q \right) | \partial_r w_{\vartheta} (r_0)|^q >0,$$
     and the contradiction. Therefore, $\gamma = r_0^\beta/\beta$ and $w= w_{r_0}$,
which completes the proof.  
\qed

\mysection{Some properties of solutions to (\ref{P})}

We now focus on time-dependent solutions to (\ref{P}) and establish some qualitative properties of non-negative and radially symmetric viscosity solutions to (\ref{P}) which are needed to analyse their large time behaviour.

\begin{prop}\label{pr3.1}
Assume that $u_0$, $p$, and $q$ fulfil (\ref{0.2}) and (\ref{0.3}). There is a unique non-negative viscosity solution $u\in \mathcal{C}(\bar{B}\times [0,\infty))$ to (\ref{P}) such that $u(x,t)=0$ for $x\in\partial B$ and $x\longmapsto u(x,t)$ is radially symmetric and belongs to $W^{1,\infty}(B)$ for all $t\ge 0$. In addition, there is a constant $A_0>0$ depending only on $p$, $q$, and $u_0$, and a decreasing function $W\in \mathcal{C}^{1}([0,\infty))$ such that 
\begin{equation}\label{3.0}
0 \le u(x,t) \le A_0 \;\;\mbox{ and }\;\; -A_0 \le \nabla u(x,t)\cdot \frac{x}{|x|} \le W(t) \,, \quad (x,t)\in \bar{B}\times [0,\infty),
\end{equation}
and $W(t)\longrightarrow 0$ as $t\to \infty$. 
\end{prop}

\proof
We first derive the expected properties on suitable approximations to (\ref{P}) which we introduce now. For $\eps\in (0,1)$, let $a_\eps\in \mathcal{C}^\infty([0,\infty))$ and $b_\eps\in \mathcal{C}^\infty([0,\infty))$ be two functions such that 
\begin{itemize}
\item $a_\eps$ is bounded and increasing and $a_\eps (\xi) := (\eps^2 + \xi)^{(p-2)/2}$ for $\xi\in [0,\eps^{-1}]$, 
\item $b_\eps$ is increasing, Lipschitz continuous, and $b_\eps (\xi) := (\eps^2 + \xi)^{q/2} - \eps^q$  for $\xi\in [0,\eps^{-1}]$.
\end{itemize}
In addition, owing to the properties (\ref{0.2}) of $u_0$, there exists a sequence $(u_{0\eps})_{\eps\in (0,1)}$ of non-negative and radially symmetric functions in $\mathcal{C}^\infty(\bar{B})$ such that
\begin{equation}\label{3.1}
\|u_{0\eps}\|_{L^\infty(B)} \le \|u_{0}\|_{L^\infty(B)} + \eps\ , \quad \|\nabla u_{0\eps}\|_{L^\infty(B)} \le 2\ \|\nabla u_{0}\|_{L^\infty(B)} ,
\end{equation}
and
$$
\lim\limits_{\eps\to 0} \| u_{0\eps}-u_0\|_{\mathcal{C}(\bar{B})} = 0 .
$$

Fix $\eps\in (0,1)$. According to the properties of $a_\eps$, $b_\eps$, and $u_{0\eps}$, it follows from \cite{LSU68} that the initial-boundary value problem
\begin{equation}\label{Peps}
  \left\{ \begin{array}{lllll}
	\partial_t u_\eps & = & \divergence(a_\eps (|\nabla u_\eps|^2) \nabla u_\eps) + b_\eps (|\nabla u_\eps|^2),	& \qquad x\in B, & \ t\in (0,\infty), \\[1mm]
	u_\eps & = & 0, & \qquad x\in\partial B, & \ t\in (0,\infty), \\[2mm]
	u_\eps |_{t=0} & = & u_{0 \eps}, & \qquad x\in B , & 
	\end{array} \right.
\end{equation} 
has a unique non-negative classical solution $u_\eps$. In addition, $x\longmapsto u_\eps(t,x)$ is radially symmetric for every $t\ge 0$ and the comparison principle entails that
\begin{equation}\label{3.2}
0 \le u_\eps(x,t) \le \|u_{0\eps}\|_{L^\infty(B)} \le \|u_{0}\|_{L^\infty(B)} + \eps, \qquad (x,t)\in \bar{B}\times [0,\infty).
\end{equation}

We next derive some estimates on the gradient of $u_\eps$ and begin with the normal trace $\partial_r u_\eps(1,t)$. Let $\mathcal{L}_\eps$ be the parabolic operator
$$
\mathcal{L}_\eps z := \partial_t z - \frac{1}{r^{N-1}} \partial_r \left( r^{N-1} a_\eps\left( |\partial_r z|^2 \right) \partial_r z \right) - b_\eps\left( |\partial_r z|^2 \right) , \quad (r,t)\in (0,1)\times (0,\infty),
$$
and fix 
\begin{equation}\label{3.2b}
A_0\in (\sqrt{3} \eps, \eps^{-1/2}) \;\;\mbox{ such that }\;\; A_0\ge 2^{1/(p-1-q)} + 2 \left( 1 + \|u_0\|_{L^\infty(B)} + \|\nabla u_0\|_{L^\infty(B)} \right).
\end{equation}
Then, thanks to the properties of $a_\eps$, $b_\eps$, and (\ref{3.2b}),  the function $\psi$ defined by $\psi(r):=A_0 (1-r)$ for $r\in [0,1]$ satisfies
\begin{eqnarray*}
\mathcal{L}_\eps \psi(r) & = & \frac{1}{r^{N-1}} \partial_r \left( r^{N-1} a_\eps\left( A_0^2 \right) A_0 \right) - b_\eps\left( A_0^2 \right) = \frac{N-1}{r}\ a_\eps\left( A_0^2 \right) A_0 - b_\eps\left( A_0^2 \right) \\
& \ge & \left( \eps^2 + A_0^2 \right)^{(p-2)/2} A_0 - \left( \eps^2 + A_0^2 \right)^{q/2} + \eps^q \\
& \ge & \left( \eps^2 + A_0^2 \right)^{(p-2)/2} \left( \sqrt{\eps^2 + A_0^2} - \eps \right) - \left( \eps^2 + A_0^2 \right)^{q/2} \\
& \ge & \left( \eps^2 + A_0^2 \right)^{(p-1)/2} \left( 1 - \frac{\eps}{\sqrt{\eps^2 + A_0^2}} \right) - \left( \eps^2 + A_0^2 \right)^{q/2} \\
& \ge & \frac{1}{2} \left( \eps^2 + A_0^2 \right)^{(p-1)/2} - \left( \eps^2 + A_0^2 \right)^{q/2} \ge 0, 
\qquad r \in (0,1].
\end{eqnarray*}
Furthermore, (\ref{3.1}), (\ref{3.2}), and (\ref{3.2b}) entail that 
$$
u_\eps\left( \frac{1}{2} , t \right) \le 1 + \|u_0\|_{L^\infty(B)} \le \frac{A_0}{2} = \psi\left( \frac{1}{2} \right)\,, \quad t\ge 0\,,
$$
and
$$
u_{0\eps}(r) = - \int_r^1 \partial_r u_{0\eps}(\varrho)\ d\varrho \le 2 \|\nabla u_0\|_\infty (1-r) \le \psi(r)\,, \quad r\in \left( \frac{1}{2} , 1 \right).
$$
Since $\mathcal{L}_\eps u_\eps=0$ in $(1/2,1)\times (0,\infty)$, the comparison principle ensures that $u_\eps(r,t)\le A_0 (1-r)$ for $(r,t)\in (1/2,1)\times (0,\infty)$. Since $u_\eps(1,t)=0$, this implies in particular that $0\le - \partial_r u_\eps(1,t)\le A_0$ for $t\ge 0$. Recalling that $u_\eps(t)$ is radially symmetric and smooth, we thus have
\begin{equation}\label{3.3}
- A_0 \le \partial_r u_\eps(1,t) \le 0=\partial_r u_\eps(0,t)\,, \quad t\ge 0\,.
\end{equation}

We next estimate the gradient of $u_\eps$ in $B$. For that purpose, we introduce the parabolic operator
\begin{eqnarray*}
\mathcal{M}_\eps z & := & \partial_t z - \partial_r \left[ \left( a_\eps\left( z^2 \right) + 2a_\eps'\left( z^2 \right)  z^2 \right) \partial_r z \right] \\
& & - \left[ \frac{N-1}{r} \left( a_\eps\left( z^2 \right) + 2a_\eps'\left( z^2 \right)  z^2 \right) + 2b_\eps'\left( z^2 \right)  z^2 \right] \partial_r z + \frac{N-1}{r^2} a_\eps\left( z^2 \right) z  
\end{eqnarray*}
for $(r,t)\in (0,1)\times (0,\infty)$ and readily deduce from (\ref{Peps}) that 
\begin{equation}\label{3.4}
\mathcal{M}_\eps \partial_r u_\eps = 0 \quad \mbox{in } (0,1)\times (0,\infty).
\end{equation}
Observe next that $\partial_r u_\eps(r,0)\ge - 2 \|\nabla u_0\|_{L^\infty(B)} \ge -A_0$ by (\ref{3.1}) and (\ref{3.2b}) and 
$$
\mathcal{M}_\eps (-A_0) = - \frac{N-1}{r^2} a_\eps\left( A_0^2 \right) A_0 \le 0\,,
$$
which, together with (\ref{3.3}), (\ref{3.4}), and the comparison principle implies that
\begin{equation}\label{3.5}
-A_0 \le \partial_r u_\eps(r,t) , \qquad (r,t)\in [0,1]\times [0,\infty) .
\end{equation}

Finally, let $W_\eps\in\mathcal{C}^1([0,\infty))$ be the solution to the ordinary differential equation
\begin{equation}\label{3.6}
\frac{dW_\eps}{dt} + (N-1)\ a_\eps\left( W_\eps^2 \right) W_\eps = 0 , \qquad W_\eps(0)=2 \|\nabla u_0\|_{L^\infty(B)} .
\end{equation}
Then $W_\eps$ is positive and decreasing, $W_\eps(0)\ge \partial_r u_\eps(r,0)$ for $r\in(0,1)$ by (\ref{3.1}), and $\mathcal{M}_\eps W_\eps\ge 0$ in $(0,1)\times (0,\infty)$ by (\ref{3.6}). Recalling (\ref{3.4}), we deduce from the comparison principle that 
\begin{equation}\label{3.7}
\partial_r u_\eps(r,t) \le W_\eps(t) , \qquad (r,t)\in [0,1]\times [0,\infty) .
\end{equation}

Finally, we argue as in \cite[Lemma~5]{GGK03} to deduce from (\ref{Peps}), (\ref{3.2}), (\ref{3.5}), and (\ref{3.7}) that there is a constant $C$ depending on $\| \nabla u_0 \|_{L^\infty (B)}$, $p$, $q$, and $N$, such that
\begin{equation}\label{3.8}
|u_\eps(x,t_1) - u_\eps(x,t_2)| \le C (|t_1-t_2| + |t_1-t_2|^{1/2})
\end{equation}
for any $x \in \bar{B}$, $t_1,t_2 \in [0,\infty)$ and $\eps \in (0,1)$. Indeed, consider $t_1 \neq t_2$ and set $\tau := 
|t_1-t_2|^{1/2} >0$ and $L:=\max\{A_0, 2 \|\nabla u_0\|_{L^\infty(B)} \}$. Since (\ref{3.5}),  (\ref{3.7}), and the Dirichlet boundary conditions imply that $|u_\eps (x,t)| \le L \, \dist(x ,\partial B)$ for $(x,t) \in \bar{B} \times [0,\infty)$, we have
\begin{equation}\label{teq2.1}
|u_\eps (x_0,t_1) - u_\eps (x_0,t_2)| \le 2 L \, \dist(x_0, \partial B) \le 2 L \tau \;\;\;\mbox{ if }\;\;\;   \dist(x_0, \partial B) \le \tau.
  \end{equation}
  If $\dist(x_0, \partial B) >\tau$ and $\eps\in (0, 1/L)$, we infer from (\ref{Peps}), the properties of $(a_\eps,b_\eps)$, and $|\nabla u_\eps | \le L$ in $B \times [0,\infty)$ that 
  \begin{eqnarray*}
    |u_\eps (x_0,t_1) - u_\eps (x_0,t_2)| &=& \frac{1}{|B| \tau^N} \left| \int\limits_{\{|x-x_0|<\tau\}} (u_\eps (x_0,t_1) 
    - u_\eps (x_0,t_2)) {\rm d}x \right| \\
    &=& \frac{1}{|B| \tau^N} \Bigg| \int\limits_{ \{ |x-x_0|<\tau \} } (u_\eps (x,t) - u_\eps (x_0,t)) 
    {\rm d}x \Big|_{t=t_1}^{t=t_2} \\ 
    & & - \int\limits_{t_1}^{t_2} \int\limits_{\{ |x-x_0|<\tau \}} \partial_t u_\eps(x,t) {\rm d}x {\rm d}t \Bigg| \\ 
    &\le & \frac{2L}{|B| \tau^N} \int\limits_{\{|x-x_0|<\tau\}} |x-x_0| {\rm d}x \\
    & & + \frac{1}{|B| \tau^N} \Bigg| \int\limits_{t_1}^{t_2} \int\limits_{\{|x-x_0|<\tau\}} \Big[ \divergence(a_\eps (|     \nabla  u_\eps|^2) \nabla u_\eps) + b_\eps (|\nabla u_\eps|^2) \Big] (x,t)  {\rm d}x {\rm d}t \Bigg| \\
    &\le& \frac{2LN}{N+1} \, \tau + \frac{1}{|B| \tau^N} \left| \int\limits_{t_1}^{t_2} \int\limits_{\{|x-x_0|<\tau\}}
    (\eps^2 + |\nabla u_\eps|^2)^{q/2}  (x,t) {\rm d}x {\rm d}t \right| \\
    & & + \frac{1}{|B| \tau^N} \left| \int\limits_{t_1}^{t_2} \int\limits_{\{|x-x_0|=\tau\}} \left[
    a_\eps (|\nabla u_\eps|^2) |\nabla u_\eps| \right] (y,t) {\rm d}S {\rm d}t \right| \\
    &\le& \frac{2LN}{N+1} \, \tau + (1+ L^2)^{q/2} |t_1-t_2|  
    + \frac{N}{\tau} (1+ L^2)^{(p-2)/2} L |t_1-t_2| \\
    & \le & \frac{2LN}{N+1}\,  \tau + (1+ L^2)^{q/2} \, \tau^2 + N (1+ L^2)^{p/2}\, \tau.
  \end{eqnarray*}
  Combining (\ref{teq2.1}) and the above estimate gives the claim (\ref{3.8}).
  
  \bigskip

We can now pass to the limit as $\eps\to 0$. Owing to (\ref{3.2}), (\ref{3.5}), (\ref{3.7}), and (\ref{3.8}), $(u_\eps)_\eps$ is bounded in, say, $C^{0,1/2}(B\times (0,\infty))$ because the uniform Lipschitz continuity in $r$ implies a uniform $C^{0,1/2}$-bound in $r$; thus $(u_\eps)_\eps$ is relatively compact in $\mathcal{C}(\bar{B}\times [0,T])$ for all $T>0$. It follows from the stability theorem \cite[Section~6]{CIL92} and the comparison principle for (\ref{P}) \cite[Theorem~2.1]{GGIS91} that $(u_\eps)_{\eps}$ converges uniformly towards the unique viscosity solution $u$ to (\ref{P}) on compact subsets of $\bar{B}\times [0,\infty)$. The properties of $u$ and the bounds listed in Proposition~\ref{pr3.1} then readily follow from this convergence, the properties of $u_\eps$, (\ref{3.2}), (\ref{3.5}), and (\ref{3.7}), the function $W$ being the solution to the ordinary differential equation 
$$
\frac{dW}{dt} + (N-1)\ |W|^{p-2} W = 0 , \qquad W(0)=2 \|\nabla u_0\|_{L^\infty(B)} .
$$
In fact, $W(t)=\left( W(0)^{2-p} + (p-2)(N-1) t \right)^{-1/(p-2)}$ for $t\ge 0$ and $W$ is obviously positive, decreasing and converges to zero as $t\to\infty$. \qed

By (\ref{3.0}), the trajectory $\{ u(t)\ : \ t\ge 0\}$ of the solution $u$ to (\ref{P}) is bounded in $L^\infty(B)$. More precise information are gathered in the next lemma.

\begin{lem}\label{le3.2}
Assume that $u_0$, $p$, and $q$ fulfil (\ref{0.2}) and (\ref{0.3}). Let $u$ be the viscosity solution to 
(\ref{P}) described in Proposition~\ref{pr3.1}. Then $t\longmapsto \|u(t)\|_{L^\infty(B)}$ is a non-increasing function and 
\begin{equation}\label{3.9}
M_\infty := \lim\limits_{t\to\infty} \|u(t)\|_{L^\infty(B)} >0 .
\end{equation}
\end{lem}

\proof
Any positive constant being obviously a supersolution to (\ref{P}), the time monotonicity of the $L^\infty(B)$-norm of $u$ readily follows from the comparison principle. Next, since $u_0\not\equiv 0$ by (\ref{0.2}), there is $x_0\in B$, $\varrho>0$, and $m>0$ such that 
$$
B_\varrho(x_0):= \{ x \in \mathbb{R}^N \ : \ |x-x_0|<\varrho \}\subset B \;\;\;\mbox{ and }\;\;\; u_0(x) \ge m \;\;\;\mbox{ for }\;\;\; x\in B_\varrho(x_0) .
$$
Introducing $v_\lambda(x):=\lambda^{(p-q)/(p-1-q)} w_0(|x-x_0|/\lambda)$ for $x\in B_\lambda(x_0)$ and $\lambda\in (0,1)$ (the function $w_0$ being defined in Theorem~\ref{theo1.1}), a simple computation shows that $v_\lambda$ is a solution to $-\Delta_p v_\lambda - |\nabla v_\lambda|^q = 0$ in $B_\lambda(x_0)$ with $v_\lambda(x)=0\le u(x,t)$ for $(x,t)\in \partial B_\lambda(x_0)\times (0,\infty)$. Furthermore, if $\lambda=\lambda_m:=\min\{ 1-|x_0| , (m\alpha/c_0)^{(p-1-q)/(p-q)}$\} , we have $v_{\lambda_m}(x)\le m \le u_0(x)$ for $x\in B_{\lambda_m}(x_0)$. The comparison principle \cite[Theorem~2.1]{GGIS91} then warrants that $u(x,t)\ge v_{\lambda_m}(x)$ for $(x,t)\in B_{\lambda_m}(x_0)\times (0,\infty)$. In particular, $\|u(t)\|_{L^\infty(B)} \ge \|v_{\lambda_m}\|_{L^\infty(B_{\lambda_m}(x_0))}$ for all $t\ge 0$, whence $M_\infty\ge \|v_{\lambda_m}\|_{L^\infty(B_{\lambda_m}(x_0))}>0$. \qed

\mysection{Convergence to steady states}

We introduce the half-relaxed limits
$$u_\ast (x) := \liminf\limits_{(s,\eps) \to (t,0)} u(x,\eps^{-1}s) , \quad x
\in \bar{B},$$
and
$$u^\ast (x) := \limsup\limits_{(s,\eps) \to (t,0)} u(x,\eps^{-1}s) , \quad x
\in \bar{B},$$
which are well-defined and do not depend on $t > 0$. 
Moreover, we infer from the stability theorem (see \cite[Lemma~6.1]{CIL92}) that
\begin{equation}\label{4.6}
  u^\ast \mbox{ is a viscosity subsolution to } F(\nabla z, D^2 z) = 0 \mbox{ in } B, 
\end{equation}
\begin{equation}\label{4.7}
  u_\ast \mbox{ is a viscosity supersolution to } F(\nabla z, D^2 z) = 0 \mbox{ in }
B. 
\end{equation}
Next we state some useful properties of the half-relaxed limits. 

\begin{lem}\label{lem4.1}
  The half-relaxed limits $u_\ast$ and $u^\ast$ enjoy the following properties:
  \begin{eqnarray}
    & & u_\ast \in W^{1,\infty} (B), \quad u^\ast \in W^{1,\infty} (B), \label{4.2} \\
    & & 0 \le u_\ast (x) \le u^\ast (x), \quad x \in \bar{B}, \label{4.3} \\
    & & u_\ast \mbox{ and } u^\ast \mbox{ are radially symmetric and non-increasing,
} \label{4.4} \\
    & & u_\ast (0) = u^\ast (0) = M_\infty := \lim\limits_{t \to \infty} \| u(t)
\|_{L^\infty (B)} >0, \label{4.5} \\
    & & u_\ast (x) = u^\ast (x) = 0 \qquad\mbox{for } x \in \partial B. \label{4.7a}
  \end{eqnarray}
\end{lem}

\proof   By (\ref{3.0}) there is $L:=\max{\{A_0,W(0)\}}>0$ such that
  \begin{equation}\label{4.1.1}
    u(x,\eps^{-1}s)  \le u(y,\eps^{-1}s) +  L |x-y| \qquad\mbox{for all } (x,y,\eps^{-1}s) \in \bar{B}\times\bar{B}\times [0,\infty)\, ,
  \end{equation}  
  from which we deduce that $u_\ast$ and $u^\ast$ are Lipschitz continuous in $B$ by taking the $\limsup$ or $\liminf$ in $\eps$
  and $s$. This proves (\ref{4.2}), while (\ref{4.3}) comes directly from the definition of $u_\ast$ and $u^\ast$ and the facts that
$u$ is non-negative, radially symmetric for any $t \ge 0$ and vanishes identically on $\partial B\times (0,\infty)$. The proof of (\ref{4.7a}) uses, in addition, the uniform Lipschitz and $C^{0,1/2}$-bounds we have for $u$ in space and time respectively.
  
 In order to prove (\ref{4.4}), we use Proposition~\ref{pr3.1}: there
  is a decreasing function $W$ such that $W(t) \to 0$ as $t \to \infty$ and
  \begin{equation}\label{4.1.2}
    u(x,t) \le u(y,t) + W(t) (|x|- |y|) \;\;\;\mbox{ for }\;\;\;  (x,y) \in \bar{B}\times\bar{B} \;\;\;\mbox{ such that }\;\;\; |x|\ge |y|.
  \end{equation}
Using this inequality with $t=\eps^{-1}s$ and taking the $\limsup$ or $\liminf$ in $\eps$ and $s$ lead to either $u_\ast (x) \le u_\ast (y)$ or
$u^\ast (x) \le u^\ast (y)$ for any $(x,y) \in \bar{B}\times\bar{B}$ such that $|x|\ge |y|$ because $W(t) \to 0$ as $t \to \infty$, hence to (\ref{4.4}).

 It remains to show (\ref{4.5}). To this end, we recall that $M_\infty$ is well-defined and positive by (\ref{3.9}) and  
  first claim that 
  \begin{equation}\label{4.1.3}
    \lim\limits_{t \to \infty} u(0,t) = M_\infty.
  \end{equation}
  Indeed, (\ref{4.1.2}) implies 
  $$u(x,t) \le u(0,t) + W(t) |x| \le u(0,t) + W(t) \le \| u(t) \|_{L^\infty (B)} + W(t), \qquad x \in B$$
  whence
  $$ \| u(t) \|_{L^\infty (B)} \le u(0,t) + W(t) \le \| u(t) \|_{L^\infty (B)} + W(t), $$
  and (\ref{4.1.3}) due to $W(t) \to 0$ as $t \to \infty$.
  
  Moreover, by the definition of the half-relaxed limits, we have $u_\ast(0) = u^\ast(0)=M_\infty$ and
  $$\|u_\ast \|_{L^\infty (B)} \le \|u^\ast \|_{L^\infty (B)} \le M_\infty.$$
  This completes the proof of (\ref{4.5}).  
\qed

Now, owing to the monotonicity and radial symmetry of $u_\ast$ and $u^\ast$, there
are $r_\ast \in [0,1]$ and $r^\ast \in
[0,1]$ such that
\begin{equation}\label{4.8}
  u_\ast (x) = M_\infty \mbox{ if } |x| \le r_\ast \mbox{ and } u_\ast (x) <
M_\infty \mbox{ if } |x| \in (r_\ast,1],
\end{equation} 
\begin{equation}\label{4.9}
  u^\ast (x) = M_\infty \mbox{ if } |x| \le r^\ast \mbox{ and } u^\ast (x) <
M_\infty \mbox{ if } |x| \in (r^\ast,1],
\end{equation}
Due to (\ref{4.3}), (\ref{4.5}), and (\ref{4.7a}), we have
\begin{equation}\label{4.10}
  0 \le r_\ast \le r^\ast < 1.
\end{equation}

Next, we show that $\Lambda u_\ast$ is a strict supersolution to the stationary
equation in a subset of $B$ for $\Lambda>1$.

\begin{lem}\label{lem4.2}
  Fix $\Lambda >1$ and $\delta \in (0, 1- r_\ast)$. Then there are $r_\delta \in
(r_\ast, r_\ast + \delta)$ and
  $\eps_{\delta, \Lambda} >0$ such that $\Lambda u_\ast$ is a viscosity
supersolution to
  $f(r, \partial_r z, \partial_r^2 z)= \eps_{\delta, \Lambda}$ in $(r_\delta,1)$.
  In addition, $\eps_{\delta, \Lambda} \to 0$ as $\Lambda \searrow 1$.
\end{lem}

\proof Fix $\delta \in (0, 1- r_\ast)$. Then, due to (\ref{4.2}), (\ref{4.4}), and
(\ref{4.8}), there is $r_\delta \in
  (r_\ast, r_\ast + \delta)$ such that $u_\ast$ is
  differentiable at $r_\delta$ and $\partial_r u_\ast (r_\delta)<0$. 
  Since $u_\ast$ is a viscosity supersolution to $f(r, \partial_r z, \partial_r^2
z)= 0$ in $(0,1)$, it is also a
  viscosity supersolution to $f_0(r, \partial_r z, \partial_r^2 z)= 0$ in $(0,1)$
and it follows from
  Lemma~\ref{lem2.1} that 
  $$\partial_r u_\ast (r) \le r^{(N-1) / (p-1)} \partial_r u_\ast (r) \le
r_\delta^{(N-1) / (p-1)} \partial_r u_\ast
    (r_\delta) =: -m_\delta <0$$
  for a.e. $r \in (r_\delta,1)$. Integrating and using the continuity of $u_\ast$ we
conclude that
  \begin{equation}\label{4.2.1}
    u_\ast (r) \le u_\ast (r_1) - m_\delta (r-r_1)
  \end{equation} 
  for all $r_1 \in [r_\delta,1]$ and $r \in [r_1,1]$. 
  
  Consider $\Lambda >1$, $\Phi \in \mathcal{C}^2 ((r_\delta,1))$ and assume that $\Lambda
u_\ast - \Phi$ has a 
  local minimum at some $r_0 \in (r_\delta,1)$. Then $u_\ast - (\Phi / \Lambda)$ has
a local minimum at
  $r_0$ and (\ref{4.7}) implies 
  \begin{eqnarray*}
  & & -\frac{1}{r_0^{N-1}} \partial_r \left( r^{N-1} \left|\partial_r \left(
\frac{\Phi}{\Lambda} \right) \right|^{p-2}
    \partial_r \left( \frac{\Phi}{\Lambda} \right) \right) (r_0)
    - \left| \partial_r \left(\frac{\Phi}{\Lambda} \right) (r_0) \right|^q
    \ge 0 , \\
    & & -\frac{1}{r_0^{N-1}} \partial_r \left( r^{N-1} \left|\partial_r \Phi \right|^{p-2}
    \partial_r \Phi \right) (r_0)
    - \Lambda^{p-1-q}\ \left| \partial_r \Phi (r_0) \right|^q
    \ge 0 .
    \end{eqnarray*}
  Thus, we have
  \begin{equation}\label{4.2.2}
    - \frac{1}{r_0^{N-1}} \partial_r \left( r^{N-1} |\partial_r \Phi|^{p-2}
\partial_r \Phi \right) (r_0)
    - | \partial_r \Phi (r_0) |^q \ge \left( \Lambda^{p-1-q} -1 \right) |\partial_r
\Phi(r_0) |^q.  
  \end{equation}
  Now, since $\Lambda u_\ast - \Phi$ has a local minimum at $r_0$, we infer from (\ref{4.2.1}) that, for $r \in [r_0, r_0+ \eta]$ with $\eta >0$ small enough,
  $$u_\ast (r_0) \le \frac{\Phi (r_0)}{\Lambda} + u_\ast (r) - \frac{\Phi(r)}{\Lambda} 
    \le \frac{\Phi (r_0)}{\Lambda} + u_\ast (r_0) - m_\delta (r-r_0) -
\frac{\Phi(r)}{\Lambda}.$$
  Hence,
  $$\frac{\Phi(r)}{\Lambda} - \frac{\Phi(r_0)}{\Lambda} \le - m_\delta (r-r_0)$$
  and thus
  $$\frac{1}{\Lambda} \partial_r \Phi(r_0) \le - m_\delta <0$$
  which implies $|\partial_r \Phi(r_0)| \ge \Lambda m_\delta$. Consequently,
(\ref{4.2.2}) becomes
  $$f(r_0, \partial_r \Phi (r_0), \partial_r^2 \Phi (r_0)) \ge \left(
\Lambda^{p-1-q} -1 \right) \Lambda^q m_\delta^q
    =: \eps_{\delta, \Lambda} >0,$$
  which ends the proof.   
\qed

We are now able to prove that the half-relaxed limits $u_\ast$ and $u^\ast$
coincide.

\begin{lem}\label{lem4.3}
  We have $u_\ast = u^\ast$ on $\bar{B}$. 
\end{lem}

\proof We fix $\Lambda > 1 > \lambda> 0$ such that $\lambda > r_\ast$ and
  $$\delta := \frac{M_\infty}{\|\nabla u_\ast \|_{L^\infty (B)}} \left( 1-
\lambda^{(p-q) / (p-1-q)} \right) \in 
  (0,\lambda - r_\ast).$$
  Defining now
  $$U(r) := \Lambda u_\ast (r) , \quad r \in [0,1], \quad\mbox{and}\quad
    V(r) := \lambda^{(p-q) / (p-1-q)} u^\ast \left( \frac{r}{\lambda} \right), \quad
r \in [0,\lambda],$$
  we obtain due to (\ref{4.10})
  \begin{equation}\label{4.3.1}
    U(r) \ge u_\ast (r) = M_\infty \ge V(r) \quad\mbox{ for } r \in [0, r_\ast].
  \end{equation}  
  Furthermore, we infer from the Lipschitz continuity of $u_\ast$ that, for $r\in (r_\ast, r_\ast + \delta]$,
  \begin{eqnarray*}
    U(r) &\ge& u_\ast (r) \ge u_\ast \left( r_\ast \right) - \| \nabla u_\ast \|
_{L^\infty (B)}
    \left| r - r_\ast \right|  \\  
    &=& M_\infty - \| \nabla u_\ast \|_{L^\infty (B)} \left| r - r_\ast \right| \ge
M_\infty - \delta \| \nabla u_\ast
    \|_{L^\infty (B)}  \\
    &\ge & \lambda^{(p-q) / (p-1-q)} M_\infty \ge V(r) .
  \end{eqnarray*}
  Recalling (\ref{4.3.1}), we have thus shown that 
  \begin{equation}\label{4.3.2}
  U(r)\ge V(r) \;\;\;\mbox{ for }\;\;\; r\in [0,r_\ast+\delta]\,.
\end{equation}
  Next, we define $I_\lambda := (r_\ast + \delta, \lambda)$. On the one hand, $V$ is a viscosity
  subsolution to
  $f(r, \partial_r z, \partial_r^2 z)= 0$ in $I_\lambda$. Indeed, take $\Phi
\in \mathcal{C}^2 (I_\lambda)$ and assume that
  $V- \Phi$ has a local maximum at $r_1 \in I_\lambda$. Then $u^\ast - \Psi$ has a
local maximum at $r_1 / \lambda$,
  where $\Psi (r) := \lambda^{-(p-q) / (p-1-q)} \Phi (\lambda r)$ for $r \in
((r_\ast + \delta) / \lambda,1)$. 
  Owing to (\ref{4.6}), we obtain
  $$f \left( \frac{r_1}{\lambda}, \partial_r \Psi \left( \frac{r_1}{\lambda}
\right), \partial_r^2 \Psi 
  \left( \frac{r_1}{\lambda} \right) \right) \le 0.$$
  Consequently,
  \begin{eqnarray*}
    0 &\ge& \lambda^{q / (p-1-q)} f \left( \frac{r_1}{\lambda}, \lambda^{-1
/(p-1-q)} \partial_r \Phi (r_1), 
    \lambda^{1- 1/ (p-1-q)} \partial_r^2 \Phi (r_1) \right)\\
    &=& - (p-1) |\partial_r \Phi (r_1)|^{p-2} \partial_r^2 \Phi (r_1) -
\frac{N-1}{r_1} |\partial_r \Phi (r_1)|^{p-2}  
    \partial_r \Phi (r_1) - |\partial_r \Phi (r_1)|^q  \\
    &=& f(r_1, \partial_r \Phi (r_1), \partial_r^2 \Phi(r_1))
  \end{eqnarray*}
  and $V$ is a viscosity subsolution to $f(r, \partial_r z, \partial_r^2
z) = 0$ in $I_\lambda$. On the other hand, it follows from  Lemma~\ref{lem4.2} that $U$ is a viscosity supersolution to $f(r,
\partial_r z, \partial_r^2 z) =
  \eps_{\delta, \Lambda}$ in $I_\lambda$ with some $\eps_{\delta,\Lambda} >0$. As
furthermore $V(r) = 0 \le U(r)$
  for $r = \lambda$ and $U(r) \ge V(r)$ for $r = r_\ast + \delta$ due to
(\ref{4.3.2}), we conclude that
  $$U(r) \ge V (r) \qquad\mbox{for } r \in [r_\ast + \delta, \lambda]$$
  by \cite[Section~5C]{CIL92}. Using (\ref{4.3.2}), we end up with
  $$\Lambda u_\ast (r) \ge \lambda^{(p-q) / (p-1-q)} u^\ast \left( \frac{r}{\lambda}
\right) \quad\mbox{ for }
    r \in [0,\lambda].$$
  Letting now $\Lambda \searrow 1$ and $\lambda \nearrow 1$, we conclude $u_\ast \ge
u^\ast$ in $[0,1]$ which, together with (\ref{4.3}),  implies $u^\ast = u_\ast$. 
\qed 

Finally, we prove Theorem~\ref{theo1.2}.

\bigskip

{\sc Proof of Theorem~\ref{theo1.2}.}

  Defining $u_\infty := u_\ast = u^\ast$ by Lemma~\ref{lem4.3}, (\ref{4.6}), (\ref{4.7}), and Lemma~\ref{lem4.1} imply that $u_\infty$ is a radially symmetric, non-increasing,
  and Lipschitz continuous viscosity solution to $F(\nabla z, D^2 z) = 0$ in $B$
satisfying $u_\infty = 0$ on $\partial B$. Moreover, $\| u_\infty \|_{L^\infty (B)}
  = M_\infty >0$ due to (\ref{4.5}). Hence, owing to Theorem~\ref{theo1.1},
there is a unique $\vartheta \in [0,1)$ such that $u_\infty = w_\vartheta$. 

  In particular, the equality $u_\ast = u^\ast$ and the definition of $u_\ast$ and
$u^\ast$ provide the uniform convergence
  of $u(t)$ towards $u^\ast = w_\vartheta$ in every compact subset of $B$ as $t \to
\infty$, see 
  \cite[Lemme~4.1]{Bl94} or \cite[Lemma~V.1.9]{BdCD97}. Combining this local
convergence with (\ref{4.2}) and (\ref{4.7a})
  gives 
  $$\lim\limits_{t \to \infty} \| u(t) - w_\vartheta\|_{\mathcal{C}(\bar{B})} =0$$
  and the claim is proved.  
\qed  

\textbf{Acknowledgements}

The authors would like to thank Olivier Ley for helpful discussions and comments.
Part of this work was done during visits of
Ph.~Lauren\c{c}ot to the Fachbereich Mathematik of the Universit\"{a}t Duisburg-Essen
and of C.~Stinner to the Institut de Math\'{e}matiques de Toulouse, Universit\'{e}
Paul Sabatier - Toulouse III. We would
like to express our gratitude for the invitation, support, and hospitality.  
%
%
%

%
%
%
%
\end{document}